\documentclass[12pt,fleqn,leqno]{article}

\author{
V.\ Oproiu \\
\and D.D.\ Poro\c sniuc\\  }
\date{}
\title{A class of K\"ahler Einstein structures on the cotangent
bundle }

\begin{document}

\maketitle
\begin{abstract}
We use some natural lifts defined on the cotangent bundle $T^*M$
of a Riemannian manifold $(M,g)$ in order to construct an almost
Hermitian structure $(G,J)$ of diagonal type. The obtained almost
complex structure $J$ on $T^*M$ is integrable if and only if the
base manifold has constant sectional curvature and the
coefficients as well as their derivatives, involved in its
definition, do fulfill a certain algebraic relation. Next one
obtains the condition that must be fulfilled in the case where the
obtained almost Hermitian structure is almost K\"ahlerian.
Combining the obtained results we get a family of K\"ahlerian
structures on $T^*M$, depending on two essential parameters. Next
we study three conditions under which the considered K\"ahlerian
structures are Einstein. In one of the obtained cases we get that
$(T^*M,G,J)$ has constant holomorphic curvature.


\end{abstract}

\section*{Introduction}

 In the study of the differential geometry of the cotangent bundle
$T^*M$ of a Riemannian manifold $(M,g)$ one uses several
Riemannian and semi-Riemannian metrics, induced from the
Riemannian metric $g$ on $M$. Among them, we may quote the metric
of Sasaki type and the complete lift of the metric $g$. On the
other hand, some notions similar to the natural lifts of $g$ to
the tangent bundle $TM$ of $M$, will induce some new Riemannian
and pseudo-Riemannian geometric structures with many nice
geometric properties. Next, one can get from $g$ some natural
almost complex structures on $T^*M$. The study of the almost
Hermitian structures induced from $g$ on $T^*M$ is an interesting
problem in the differential geometry of the cotangent bundle.

In the present paper we study some classes of natural  K\"ahler
Einstein structures $(G,J)$, of diagonal type induced on $T^*M$
from the Riemannian metric $g$. They are obtained in a manner
quite similar to that used in \cite{Oproiu5} (see also
\cite{OprPor}) but the parametrization is a bit different. Namely,
we adapt the situation presented in \cite{Oproiu6} to the case of
the cotangent bundle, restricting ourselves to the case of the
lifts of diagonal type. In fact we do not consider the most
general situation due to the hard computations that must be done.
However, in principle, the results obtained in the case of the
general natural almost Hermitian structures on $T^*M$ do not
differ too much from that obtained in the case of the natural
almost Hermitian structures of diagonal type. We consider the case
where the vertical and horizontal distributions are orthogonal to
each other but the dot products induced on them from $G$ are not
isomorphic (isometric). The family of the natural almost complex
structures $J$ on $T^*M$ that interchange the vertical and
horizontal distributions depends on two essential parameters $a_1,
b_1$. These parameters are smooth real functions depending on the
energy density $t$ on the cotangent bundle. From the integrability
condition for $J$ it follows that the base manifold $M$ must have
constant curvature $c$ and the second parameter $b_1$ must be
expressed as a rational function depending on the first parameter
$a_1$ and its derivative. Of course, in the obtained formula there
are involved too the constant $c$ and the energy density $t$.

A natural Riemannian metric $G$ of diagonal type on $T^*M$ is
defined by four parameters $c_1,c_2, d_1, d_2$ which are smooth
functions of $t$. From the condition for $G$ to be Hermitian with
respect to $J$ we get two sets of proportionality relations, from
which one obtains the parameters $c_1,c_2, d_1, d_2$ as functions
depending on two new parameters $\lambda, \mu$ and the parameters
$a_1, b_1 $ involved in the expression of $J$. In the case where
the fundamental $2$-form $\phi$, associated to the almost complex
structure $(G,J)$ is closed, one finds that $\mu =
\lambda^\prime$. If the the integrability condition for $J$ is
fulfilled, we get a K\"ahlerian structure on $T^*M$ and this
structure depends on two essential parameters $a_1$ and $\lambda$.

In the case where the considered K\"ahlerian structure is Einstein
we get several situations in which  the parameters $a_1, \lambda $
are related by some algebraic relations. We have a general case,
when $(T^*M,G,J)$ has constant holomorphic curvature. In other two
cases one obtains some simpler expressions for the components of
the curvature tensor field on $T^*M$ and, of course, we have some
singularities. These cases will be discussed in some forthcoming
papers.

The manifolds, tensor fields and geometric objects we consider in
this paper, are assumed to be differentiable of class $C^{\infty}$
(i.e. smooth). We use the computations in local coordinates but
many results from this paper may be expressed in an invariant
form. The well known summation convention is used throughout this
paper, the range for the indices $h,i,j,k,l,r,s$ being
always${\{}1,...,n{\}}$ (see \cite{GheOpr}, \cite{OprPap1},
\cite{OprPap2}). We shall denote by ${\Gamma}(T^*M)$ the module of
smooth vector fields on $T^*M$.

\section{Natural almost complex structures of diagonal type on
$T^*M$}

Let $(M,g)$ be a smooth $n$-dimensional Riemannian manifold and
denote its cotangent bundle by $\pi :T^*M\longrightarrow M$.
Recall that there is a structure of a $2n$-dimensional smooth
manifold on $T^*M$, induced from the structure of smooth
$n$-dimensional manifold  of $M$. From every local chart
$(U,\varphi )=(U,x^1,\dots ,x^n)$  on $M$, it is induced a local
chart $(\pi^{-1}(U),\Phi )=(\pi^{-1}(U),q^1,\dots , q^n,$
$p_1,\dots ,p_n)$, on $T^*M$, as follows. For a cotangent vector
$p\in \pi^{-1}(U)\subset T^*M$, the first $n$ local coordinates
$q^1,\dots ,q^n$ are  the local coordinates $x^1,\dots ,x^n$ of
its base point $x=\pi (p)$ in the local chart $(U,\varphi )$ (in
fact we have $q^i=\pi ^* x^i=x^i\circ \pi, \ i=1,\dots n)$. The
last $n$ local coordinates $p_1,\dots ,p_n$ of $p\in \pi^{-1}(U)$
are the vector space coordinates of $p$ with respect to the
natural basis $(dx^1_{\pi(p)},\dots , dx^n_{\pi(p)})$, defined by
the local chart $(U,\varphi )$,\ i.e. $p=p_idx^i_{\pi(p)}$. Due to
this special structure of differentiable manifold for $T^*M$ it is
possible to introduce the concept of $M$-tensor field on it. An
$M$-tensor field of type $(r,s)$ on $T^*M$ is defined by sets of
$n^{r+s}$ components (functions depending on $q^i$ and $p_i$),
with $r$ upper indices and $s$ lower indices, assigned to induced
local charts $(\pi^{-1}(U),\Phi )$ on $T^*M$, such that the local
coordinate change rule is that of the local coordinate components
of a tensor field of type $(r,s)$ on the base manifold $M$, when a
change of local charts on $M$ (and hence on $T^*M$) is performed
(see \cite{Mok} for further details in the case of the tangent
bundle); e.g., the components $p_i,\ i=1,\dots ,n$, corresponding
to the last $n$ local coordinates of a cotangent vector $p$,
assigned to an induced local chart $(\pi^{-1}(U), \Phi )$ define
an $M$-tensor field of type $(0,1)$ on $T^*M$. A usual tensor
field of type $(r,s)$ on $M$ may be thought of as an $M$-tensor
field of type $(r,s)$ on $T^*M$. If the considered tensor field on
$M$ is covariant only, the corresponding $M$-tensor field on
$T^*M$ may be identified with the induced (pullback by $\pi $)
tensor field on $T^*M$. Some useful $M$-tensor fields on $T^*M$
may be obtained as follows. Let $u:[0,\infty ) \longrightarrow
{\bf R}$ be a smooth function and let
$\|p\|^2=g^{-1}_{\pi(p)}(p,p)$ be the square of the norm of the
cotangent vector $p\in \pi^{-1}(U)$ ($g^{-1}$ is the tensor field
of type (2,0) having as components the entries $g^{ij}(x)$ of the
inverse of the matrix $(g_{ij}(x))$ defined by the components of
$g$ in the local chart $(U,\varphi )$). If $\delta ^i_j$ are the
Kronecker symbols (in fact, they are the local coordinate
components of the identity tensor field $I$ on $M$), then the
components $u(\|p\|^2)\delta ^i_j$ define an $M$-tensor field of
type $(1,1)$ on $T^*M$. Similarly, if $g_{ij}(x)$ are the local
coordinate components of the metric tensor field $g$ on $M$ in the
local chart $(U,\varphi )$, then the components $u(\|p\|^2)
g_{ij}(\pi(p))$ define a symmetric
 $M$-tensor field of type $(0,2)$ on $T^*M$. The components
 $g^{0i}=p_hg^{hi}$, as well as $u(\|p\|^2)g^{0i}$ define $M$-tensor
 fields of type $(1,0)$ on $T^*M$. Of course,  all the components
 considered above are in the induced local chart $(\pi^{-1}(U),\Phi)$.

We shall use the horizontal distribution $HT^*M$, defined by the
Levi Civita connection $\dot \nabla $ of $g$, in order to define
some first order natural lifts to $T^*M$ of the Riemannian metric
$g$ on $M$. Denote by $VT^*M= {\rm Ker}\ \pi _*\subset TT^*M$ the
vertical distribution on $T^*M$. Then we have the direct sum
decomposition

\begin{equation}
TT^*M=VT^*M\oplus HT^*M.
\end{equation}

If $(\pi^{-1}(U),\Phi)=(\pi^{-1}(U),q^1,\dots ,q^n,p_1,\dots
,p_n)$ is a local chart on $T^*M$, induced from the local chart
$(U,\varphi )= (U,x^1,\dots ,x^n)$, the local vector fields
$\frac{\partial}{\partial p_1}, \dots , \frac{\partial}{\partial
p_n}$ on $\pi^{-1}(U)$ define a local frame for $VT^*M$ over $\pi
^{-1}(U)$ and the local vector fields $\frac{\delta}{\delta
q^1},\dots ,\frac{\delta}{\delta q^n}$ define a local frame for
$HT^*M$ over $\pi^{-1}(U)$, where
$$
\frac{\delta}{\delta q^i}=\frac{\partial}{\partial
q^i}+\Gamma^0_{ih} \frac{\partial}{\partial p_h},\ \ \ \Gamma
^0_{ih}=p_k\Gamma ^k_{ih}
 $$
and $\Gamma ^k_{ih}(\pi(p))$ are the Christoffel symbols of $g$.

The set of vector fields $(\frac{\partial}{\partial p_1},\dots
,\frac{\partial}{\partial p_n}, \frac{\delta}{\delta q^1},\dots
,\frac{\delta}{\delta q^n})$ defines a local frame on $T^*M$,
adapted to the direct sum decomposition (1). Remark that
$$
\frac{\partial}{\partial p_i}=(dx^i)^V,\ \ \frac{\delta}{\delta
q^i}=(\frac{\partial}{\partial x^i})^H,
$$
where $\theta^V$ denotes the vertical lift to $T^*M$ of the
$1$-form $\theta$ on $M$ and $X^H$ denotes the horizontal lift to
$T^*M$ of the vector field $X$ on $M$.

Now we shall present the following auxiliary result. \vskip5mm

{\bf Lemma 1}. \it If $n>1$ and $u,v$ are smooth functions on
$T^*M$ such that
$$
u g_{ij}+v p_ip_j=0,\ p\in \pi^{-1}(U)
$$
on the domain of any induced local chart on $T^*M$, then $u=0,\
v=0$.\rm \vskip 0.5cm

The proof is obtained easily by transvecting the given relation
with components $g^{ij}$ of the tensor field $g^{-1}$ and $g^{0j}$
(Recall that the functions $g^{ij}(x)$ are the components of the
inverse of the matrix $(g_{ij}(x))$, associated to $g$ in the
local chart $(U,\varphi )$ on $M$).\vskip5mm

\bf Remark. \rm From the relations of the type
$$
u g^{ij}+v g^{0i}g^{0j}=0,\ p\in \pi^{-1}(U),
$$
$$
u\delta ^i_j+vg^{0i} p_j=0,\ p\in \pi^{-1}(U),
$$
it is obtained, in a similar way, $u=v=0$. We have used the
notation $g^{0i}=p_hg^{hi}$. \vskip5mm

Since we work in a fixed local chart $(U,\varphi )$ on $M$ and in
the corresponding induced local chart $(\pi^{-1}(U),\Phi )$ on
$T^*M$, we shall use the following simpler notations
$$
\frac{\partial}{\partial p_i}=\partial ^i,\ \ \frac{\delta}{\delta
q^i}= \delta _i.
$$

Denote by
\begin{equation}
t=\frac{1}{2}\|p\|^2=\frac{1}{2}g^{-1}_{\pi(p)}(p,p)=\frac{1}{2}g^{ik}(x)p_ip_k,
\ \ \ p\in \pi^{-1}(U)
\end{equation}
the energy density defined by $g$ in the cotangent vector $p$. We
have $t\in [0,\infty)$ for all $p\in T^*M$. For a vector field $X$
on $M$ we shall denote by $g_X$ the $1$-form on $M$ defined by
$g_X(Y)=g(X,Y)$, for all vector fields $Y$ on $M$. For a $1$-form
$\theta$ on $M$, we shall denote by $\theta^\sharp=g^{-1}_\theta$
the vector field on $M$ defined by the usual musical isomorphism,
i.e. $g(\theta^\sharp,Y)=\theta (Y)$, for all vector fields $Y$ on
$M$. Remark that, for $p\in T^*M$, we can consider the vector
$p^\sharp$, tangent to $M$ in $\pi(p)$. Consider the real valued
smooth functions $a_1,a_2,b_1,b_2$ defined on $[0,\infty)\subset
{\bf R}$ and define a diagonal natural almost complex structure
$J$ on $T^*M$, by using these coefficients and the Riemannian
metric $g$

\begin{equation}
\left\{
\begin{array}{l}
JX^H_p=a_1(t) (g_X)_p^V+b_1(t)p(X)p_p^V,
\\ \mbox{ }  \\
J\theta^V_p=-a_2(t)(\theta^\sharp)_p^H- b_2(t)g^{-1}_{\pi(p)}
(p,\theta)(p^\sharp)_p^H.
\end{array}
\right.
\end{equation}

We should remark that the vector $p_p^V$ defines the Liouville
vector field on $T^*M$ and $(p^\sharp)_p^H$ defines a similar
$HT^*M$-valued vector field.

The expression of $J$ in adapted local frames is given by

$$
J{\delta _i}=a_1(t)g_{ij}{\partial ^j}+ b_1(t)p_ip_j{\partial^j},
$$
$$
J{\partial ^i}=-a_2(t)g^{ij}{\delta _j}-b_2(t)g^{0i}g^{0j}\delta
_j.
$$

Remark that one can consider the  case of the general natural
tensor fields $J$ on $T^*M$. In this case we have another four
coefficients $a_3,\ b_3,\ a_4,\ b_4$ and the computations involved
in the study of the corresponding almost complex structure $J$ on
$T^*M$ become really complicate (see \cite{Oproiu6},
\cite{Oproiu2}). In fact, the tensor fields of this type define
the most general natural lift of type $(1,1)$ of the metric $g$.
\vskip5mm

{\bf Proposition 2}. \it The operator $J$ defines an almost
complex structure on $T^*M$ if and only if

\begin{equation}
a_1a_2=1\ ,\ \ \ (a_1+2tb_1)(a_2+2tb_2)=1.
\end{equation}
\vskip2mm

 Proof. \rm The relations are obtained easily from the property $J^2=-I$ of
$J$ and Lemma 1.

From the relations obtained in Lemma 1 we can get the explicit
expressions of the parameters $a_2,b_2$
\begin{equation}
a_2 =\frac{1}{a_1}, \ \ b_2=-\frac{b_1}{a_1(a_1+2tb_1)}
\end{equation}

The obtained almost complex structures defined by the tensor field
$J$ on $T^*M$ are called \it natural almost complex structures of
diagonal type, \rm defined by the Riemannian metric $g$, by using
the essential parameters $a_1,b_1$. We use the word diagonal for
these almost complex structures, since the $2n\times 2n$-matrix
associated to $J$, with respect to the adapted local frame
$(\frac{\delta}{\delta q^1},\dots ,\frac{\delta}{\delta
q^n},\frac{\partial}{\partial p_1},\dots ,\frac{\partial}{\partial
p_n})$ has two $n\times n$-blocks on the second diagonal
\begin{displaymath}
J= \left(
\begin{array}{cc}
0 & -a_2g^{ij}-b_2g^{0i}g^{0j} \\
a_1g_{ij}+b_1p_ip_j & 0
\end{array}
\right).
\end{displaymath}

\bf Remark. \rm From the conditions (4) we have that the
coefficients $a_1,a_2, a_1+2tb_1,a_2+2tb_2$ cannot vanish and have
the same sign. We
 assume that $a_1>0,\ a_2>0,\ a_1+2tb_1>0,\ a_2+2tb_2>0$ for all $t\geq 0$.

Now we shall study the integrability of the almost complex
structure defined by $J$ on $T^*M$. To do this we need the
following well known formulas for the brackets of the vector
fields $\partial^i=\frac{\partial}{\partial p_i},\delta_i=
\frac{\delta}{\delta q^i},~ i=1,...,n$
\begin{equation}
[\partial^i,\partial^j]=0;~~~[\partial^i,\delta_j]=\Gamma^i_{jk}\partial^k;~~~
[\delta_i,\delta_j] =R^0_{kij}\partial^k,
\end{equation}
where $\Gamma^i_{jk}$ are the Christoffel symbols defined by the
Levi Civita connection $\dot\nabla$,  $R^0_{kij}=p_hR^h_{kij}$ and
$R^h_{kij}$ are the local coordinate components of the curvature
tensor field of $\dot \nabla$ on $M$.\vskip2mm

{\bf Theorem 3. } {\it The almost complex structure $J$ on $T^*M$
is integrable if and only if $(M,g)$ has  constant sectional
curvature $c$ and the function $b_1$ is given by
\begin{equation}
b_1=\frac{a_1a_1^\prime- c}{a_1-2ta_1^\prime}.
\end{equation}} \vskip2mm

Of course we have to study the conditions under which $a_1,b_1$
fulfill the conditions $a_1>0,\ a_1+2tb_1 =
\frac{a_1^2-2ct}{a_1-2ta_1^\prime}>0,\ \forall t\geq 0$. \vskip2mm

{\it Proof. } We shall study the vanishing of the Nijenhuis tensor
field $N=N_J$ of $J$, defined by
$$
N(X,Y)=[JX,JY]-J[JX,Y]-J[X,JY]-[X,Y],\ \ \forall X,Y \in \Gamma
(T^*M).
$$
We have $\delta_kt =0,\ \partial^k t = g^{0k}$ and, after a
straightforward but quite long computation, we get
$$
N(\delta_i,\delta_j)= \{(a_1a_1^\prime+2t a_1^\prime b_1-a_1b_1)
(p_ig_{jk}-p_jg_{ik})-R^0_{kij}\}\partial^k
$$

Remark that the coefficient of $\delta_k$ in the expression of
$N(\delta_i,\delta_j)$ becomes $0$, due to the usual properties of
the Levi Civita connection $\dot \nabla$.

 From the condition $N(\delta_i,\delta_j)=0$ we get
$$
R^0_{kij}= (a_1a_1^\prime+2t a_1^\prime b_1-a_1b_1)
(p_ig_{jk}-p_jg_{ik}).
$$

Differentiating with respect to $p_h$, taking $p=0$ and using
Schur theorem, it follows that the curvature tensor field of $\dot
\nabla$ (in the case where $M$  is connected and ${\rm dim}\ M>2$)
must have the expression
$$
R^h_{kij}=c({\delta}^h_ig_{kj}-{\delta}^h_jg_{ki}),
$$
where $c$ is a constant. Then we obtain the expression (7) of
$b_1$.

Next it follows by a straightforward computation that
$N(\partial^i,\delta_j)=0,$\break $N(\partial^i,\partial^j)=~0,$
whenever $N(\delta_i,\delta_j)=~0$.

Hence the condition $N=0$ implies that $(M,g)$ must have constant
sectional curvature $c$, and $b_1$ must be given by (7).
Conversely, if $(M,g)$ has constant curvature $c$ and $b_1$ is
given by (7), it follows in a straightforward way that $N=0$.
\vskip2mm

\bf Remark. \rm In the case where $a_1^2-2ct=0$, we have
$a_1a_1^\prime- c=0,\ a_1-2ta_1^\prime=0$ too. So, this case must
be thought of as a singular case and should be considered
separately.

 \section{Natural almost Hermitian structures on $T^*M$}
 Consider the following symmetric $M-$tensor
fields on $T^*M$, defined by the components
\begin{equation}
G^{(1)}_{ij}=c_1 g_{ij}+d_1p_ip_j,\ \ \
G_{(2)}^{ij}=c_2g^{ij}+d_2g^{0i} g^{0j},
\end{equation}
where $c_1,c_2,d_1,d_2$ are smooth functions depending on the
energy density $t\in [0,\infty)$.

Obviously, $G^{(1)}$ is of type $(0,2)$ and $G_{(2)}$ is of type
$(2,0)$. We shall assume that the matrices defined by $G^{(1)}$
and $G_{(2)}$ are positive definite. This happens iff
$c_1>0,c_2>0,c_1+2td_1>0, c_2+2td_2>0$. Then the following
Riemannian metric may be considered on $T^*M$
\begin{equation}
G=G^{(1)}_{ij}dq^idq^j+G_{(2)}^{ij}Dp_iDp_j,
\end{equation}
where $Dp_i=dp_i-\Gamma^0_{ij}dq^j$ is the absolute (covariant)
differential of $p_i$ with respect to the Levi Civita connection
$\dot\nabla$ of $g$ (recall that
$\Gamma^0_{ij}=p_h\Gamma^h_{ij}$). Equivalently, we have
$$
G(\delta_i,\delta_j)=G^{(1)}_{ij},~~~G(\partial^i
,\partial^j)=G_{(2)}^{ij},~~ G(\partial^i,\delta_j)=~
G(\delta_j,\partial^i)=0.
$$
Remark that $HT^*M,~VT^*M$ are orthogonal to each other with
respect to $G$, but the Riemannian metrics induced from $G$ on
$HT^*M,~VT^*M$ are not the same, so the considered metric $G$ on
$T^*M$ is not a metric of Sasaki type. The $2n\times 2n$-matrix
associated to $G$, with respect to the adapted local frame
$(\frac{\delta}{\delta q^1},\dots ,\frac{\delta}{\delta
q^n},\frac{\partial}{\partial p_1},\dots ,\frac{\partial}{\partial
p_n})$ has two $n\times n$-blocks on the first diagonal
\begin{displaymath}
G= \left(
\begin{array}{cc}
G^{(1)}_{ij} & 0  \\
0 & G_{(2)}^{ij}
\end{array}
\right).
\end{displaymath}

The Riemannian metric $G$ is called a \it natural lift of diagonal
type \rm of $g$. Remark also that the system of 1-forms
$(dq^1,...,dq^n,Dp_1,...,Dp_n)$ defines a local frame on
$T^{*}T^*M$, dual to the local frame $(\delta_1 ,...,\delta_n,~
\partial^1,...,\partial^n)$ adapted to the
direct sum decomposition (1).

We shall consider another two $M$-tensor fields $H_{(1)}, \
H^{(2)}$ on $T^*M$, defined by the components
$$
H_{(1)}^{jk}=\frac{1}{c_1}g^{jk}-\frac{d_1}{c_1(c_1+2td_1)}g^{0j}g^{0k},
$$
$$
H^{(2)}_{jk}=\frac{1}{c_2}g_{jk}-\frac{d_2}{c_2(c_2+2td_2)}p_jp_k.
$$

The components $H_{(1)}^{jk}$ define an $M$-tensor field of type
$(2,0)$ and the components $H^{(2)}_{jk}$ define an $M$-tensor
field of type $(0,2)$. Moreover, the matrices associated to
$H_{(1)}, \ H^{(2)}$ are the inverses of the matrices associated
to  $G^{(1)}$ and  $G_{(2)}$, respectively, i.e. we have
$$
G^{(1)}_{ij}H_{(1)}^{jk} = \delta_i^k,\ \ G_{(2)}^{ij}H^{(2)}_{jk}
= \delta^i_k.
$$

Now, we shall be interested in the conditions under which the
metric $G$ is almost Hermitian with respect to the almost complex
structure $J$, considered in the previous section, i.e.
$$
G(JX,JY)=G(X,Y),
$$
for all vector fields $X,Y$ on $T^*M$.

Considering the coefficients of $g_{ij}, g^{ij}$ in the conditions
\begin{equation}
\left\{
\begin{array}{l}
G(J\delta_i,J\delta_j)= G(\delta_i,\delta_j),
\\ \mbox{ } \\
G(J\partial^i,J\partial^j)= G(\partial^i,\partial^j),
\end{array}
\right.
\end{equation}
we can express the parameters $c_1,c_2$  with the help of the
parameters $a_1, a_2$ and a proportionality factor $\lambda =
\lambda (t)$
\begin{equation}
c_1 = \lambda a_1,\ \ c_2=\lambda a_2,
\end{equation}
where the coefficients $a_1,a_2$ are related by (4). Since we made
the assumption $a_1>0,\ a_2>0$, it follows $\lambda>0$.

Next, considering the coefficients of $p_ip_j,\ g^{0i}g^{0j}$ in
the relations (10), we can express the parameters
$c_1+2td_1,c_2+2td_2$ with the help of the parameters $a_1+2tb_1,
a_2+2tb_2$ and a new parameter $\lambda+2t\mu $

\begin{equation}
\left\{
\begin{array}{l}
c_1+2td_1=(\lambda+2t\mu )(a_1+2tb_1),
\\ \mbox{ } \\
c_2+2td_2=(\lambda+2t\mu )(a_2+2tb_2).
\end{array}
\right.
\end{equation}
Remark that $\lambda+2t\mu =\lambda (t)+2t \mu (t)$ is a positive
smooth function of $t\in [0,\infty)$. It was much more convenient
to consider the proportionality factor in such a form in the
expression of the parameters $c_1+2td_1,c_2+2td_2$.  Of course, we
can obtain easily from (12) the explicit expressions of the
coefficients $d_1,d_2$
\begin{equation}
\left\{
\begin{array}{l}
d_1=\lambda b_1+\mu (a_1+2tb_1),
\\ \mbox{ } \\
d_2=\lambda b_2+\mu (a_2+2tb_2).
\end{array}
\right.
\end{equation}

Hence we may state

{\bf Theorem 4.} \it Let $J$ be the natural, almost complex
structure of diagonal type on $T^*M$, given by (3), where the
coefficients $a_1,a_2,b_1, b_2$ are related by (4). The family of
the natural Riemannian metrics $G$ on $T^*M$, of diagonal type,
such that $(T^*M,G,J)$ is an almost Hermitian manifold, is given
by (9) where the coefficients $c_1,c_2$ are related to $a_1,a_2$
by (11), and $c_1+2td_1,c_2+2td_2$ are related to
$a_1+2tb_1,a_2+2tb_2$ by (12), the proportionality coefficients
being $\lambda >0 $ and $\lambda +2t\mu>0$.

\bf Remark. \rm A result of the same kind can be obtained in the
case of the natural almost Hermitian structures of general type on
$T^*M$ (see \cite{Oproiu6}).

 Consider now the two-form $\phi $ defined by the almost
Hermitian structure $(G,J)$ on $T^*M$
$$
\phi (X,Y)=G(X,JY),
$$
for all vector fields $X,Y$ on $T^*M$. \vskip2mm

\bf Proposition 5. \it The expression of the $2$-form $\phi $ in a
local adapted frame $(\partial^1,\dots , \partial^n,
\delta_1,\dots ,\delta_n)$ on $T^*M$, is given by
$$
\phi (\partial^i,\partial^j)=0,\ \phi(\delta_i,\delta_j)=0,\
\phi(\partial^i,\delta_j)= \lambda \delta^i_j+\mu g^{0i}p_j,
$$
or, equivalently
\begin{equation}
\phi =(\lambda \delta^i_j+\mu g^{0i}p_j)Dp_i\wedge dq^j,
\end{equation}
where  $Dp_i=dp_i-\Gamma^0_{ih}dq^h$ is the absolute differential
of $p_i$. \rm \vskip2mm

The proof is obtained by using the definition of $\phi$ and
computing the values $\phi (\partial^i,\partial^j),
\phi(\delta_i,\delta_j), \phi(\partial^i,\delta_j)$. \vskip2mm

\bf Theorem 6. \it The almost Hermitian structure $(G,J)$ on
$T^*M$ is almost K\"ahlerian if and only if
$$
\mu=\lambda ^\prime .
$$

Proof. \rm We shall study the vanishing of the exterior
differential $d\phi$. The expressions of $d\lambda,\ d\mu, \
dg^{0i}$ and $dDp_i$ are obtained in a straightforward way, by
using the property $\dot \nabla_k g_{ij}=0$ (hence $\dot \nabla_k
g^{ij}=0$)
$$
d\lambda = \lambda ^\prime g^{0i}Dp_i,\ d\mu =\mu ^\prime
g^{0i}Dp_i,\ dg^{0i}=g^{ik}Dp_k-g^{0h}\Gamma ^i_{hk}dq^k,
$$
$$
dDp_i=-\frac{1}{2}R^0_{ikl}dq^k\wedge dq^l+ \Gamma
^l_{ik}dq^k\wedge Dp_l.
$$
Then we have
$$
d\phi =(d\lambda \delta^i_j+d\mu g^{0i}p_j+ \mu dg^{0i}p_j+\mu
g^{0i}dp_j)\wedge Dp_i\wedge dq^j+
$$
$$
+(\lambda \delta^i_j+\mu g^{0i}p_j)dDp_i\wedge dq^j.
$$
By replacing the expressions of $d\lambda , d\mu , dg^{0i}$ and
$dDp_i$ then using, again, the property $\dot\nabla _kg_{ij}=0$,
doing some algebraic computations with the exterior products, then
using the well known symmetry properties of $g_{ij}, \Gamma
^h_{ij},$ and of the Riemann-Christoffel tensor field, as well as
the Bianchi identities, it follows that
$$
d\phi =\frac{1}{2}(\lambda ^\prime -\mu)g^{0h}Dp_h\wedge
Dp_i\wedge dq^i.
$$
Therefore we have $d\phi =0$ if and only if $\mu =\lambda ^\prime
$.

\bf Theorem 7. \it The almost Hermitian structure $(G,J)$ on
$T^*M$ is K\"ahlerian if and only if the base manifold $M$ has
constant sectional curvature, the parameter $b_1$ is given by (7)
and $\mu =\lambda ^\prime$.

Proof. \rm The family of natural almost Hermitian structures
$(G,J)$ of diagonal type on $T^*M$ depends on four essential
coefficients $a_1, b_1,\lambda ,\mu$. According to the result of
theorem 3, the integrability of $J$ is equivalent to the property
of $M$ to have constant sectional curvature $c$ and the condition
for $b_1$ to be given by (7). Then, from theorem 6, we get the
$(G,J)$ is almost K\"ahlerian if and only if $\mu=\lambda^\prime$.
Combining these two results one obtains the result of our theorem.

\bf Remark. \rm A natural K\"ahlerian structure $(G,J)$ of
diagonal type on $T^*M$ is defined by two essential coefficients
$a_1,\lambda$. Using $(7)$, these coefficients must satisfy the
supplementary conditions $a_1>0, \
a_1+2tb_1=\frac{a_1^2-2ct}{a_1-2ta_1^\prime}>0,\ \lambda
>0,\ \lambda+2t\lambda^\prime>0$.

\section{The Levi Civita connection of the metric $G$ and its
curvature tensor field}

Recall that the Levi Civita connection $\dot\nabla$ on a
Riemannian manifold $(M,g)$ is determined by the conditions
$$
\dot\nabla g=0,~~~~~\dot T=0,
$$
where $\dot T$ is its torsion tenor field. The explicit expression
of this connection is obtained from the formula
$$
2g({\dot\nabla}_XY,Z)=X(g(Y,Z))+Y(g(X,Z))-Z(g(X,Y))+
$$
$$
+g([X,Y],Z)-g([X,Z],Y)-g([Y,Z],X), ~~~~~~\forall
X,Y,Z{\in}{\Gamma}(M).
$$

We shall use this formula in order to obtain the expression of the
Levi Civita connection ${\nabla}$ of $G$ on $T^*M$. The final
result can be stated as follows \vskip2mm

{\bf Theorem 8.} {\it The Levi Civita connection ${\nabla}$ of $G$
has the following expression in the local adapted frame
$(\partial^1,...,\partial^n,~ \delta_1,..., \delta_n)$
$$
\nabla_{\partial^i}\partial^j =Q^{ij}_h\partial^h,\ \ \
\nabla_{\delta_i}\partial^j=-\Gamma^j_{ih}\partial^h+P^{hj}_i\delta_h,
$$
$$
\nabla_{\partial^i}\delta_j=P^{hi}_j\delta_h,\ \ \
\nabla_{\delta_i}\delta_j=\Gamma^h_{ij}\delta_h+S_{hij}\partial^h,
$$
where $Q^{ij}_h, P^{hi}_j, S_{hij}$ are $M$-tensor fields on
$T^*M$, defined by
$$
Q^{ij}_h = \frac{1}{2}H^{(2)}_{hk}(\partial^iG_{(2)}^{jk}+
\partial^jG_{(2)}^{ik} -\partial^kG_{(2)}^{ij}),
$$
$$
P^{hi}_j=\frac{1}{2}H_{(1)}^{hk}(\partial^iG^{(1)}_{jk}-G_{(2)}^{il}R^0_{ljk}),
$$
$$
S_{hij}=-\frac{1}{2}H^{(2)}_{hk}\partial^kG^{(1)}_{ij}+\frac{1}{2}R^0_{hij}.
$$ \rm

Replacing the expressions of the involved $M$-tensor fields and
assuming that the base manifold $(M,g)$ has constant sectional
curvature, one obtains

$$
Q^{ij}_h = -\frac{c_2^\prime - 2d_2}{2(c_2 + 2d_2t)}g^{ij}p_h +
 \frac{c_2^\prime}{2c_2}(\delta_h^jg^{0i} +\delta_h^ig^{0j})+
 \frac{-2c_2^\prime d_2 + c_2d_2^\prime}{2c_2(c_2 +
 2d_2t)}p_hg^{0i}g^{0j},
$$
$$
P^{hi}_j=-\frac{cc_2 - d_1}{2c_1}g^{hi}p_j +
  \frac{cc_2 + d_1}{2(c_1 + 2d_1t)}\delta^i_jg^{0h} +
 \frac{c_1^\prime}{2c_1}\delta^h_j g^{0i} +
 $$
 $$
+ \frac{-c_1^\prime d_1 + cc_2d_1 - d_1^2 + c_1d_1^\prime}
 {2c_1(c_1 + 2d_1t)}p_jg^{0h}g^{0i},
$$
$$
S_{hij}= \frac{-c_1^\prime}{2(c_2 + 2d_2t)}g_{ij}p_h +
  \frac{cc_2 - d_1}{2c_2}g_{hj}p_i -
  \frac{cc_2 + d_1}{2c_2}g_{hi}p_j -
  $$
  $$
 - \frac{c_2d_1^\prime - 2d_1d_2}{2c_2(c_2 + 2d_2t)}p_hp_ip_j.
$$

In the case of a K\"ahler structure on $T^*M$, the final
expressions of these $M$-tensor fields can be obtained by doing
the necessary replacements. However, the final expressions are
quite complicate but they may be obtained quite automatically by
using the Mathematica package RICCI for doing tensor computations
(see \cite{Lee}}).

Now we shall indicate the obtaining of the components of the
curvature tensor field of the connection $\nabla$.

The curvature tensor $K$ field of the connection $\nabla $ is
obtained from the well known formula
$$
K(X,Y)Z=\nabla_X\nabla_YZ-\nabla_Y\nabla_XZ-\nabla_{[X,Y]}Z.
$$

The components of $K$ with respect to the adapted local frame
$(\partial^1,\dots ,$ $\partial^n,\delta_1,\dots ,\delta_n)$ can
be expressed easily
$$
K(\partial^i,\partial^j)\partial^k=PPP^{ijk}_h\partial^h
=(\partial^iQ^{jk}_h-\partial^jQ^{ik}_h+Q^{jk}_lQ^{il}_h-
Q^{ik}_lQ^{jl}_h)\partial^h,
$$
$$
K(\partial^i,\partial^j)\delta_k=PPQ^{ijh}_k\delta_h
=(\partial^iP^{hj}_k-\partial^jP^{hi}_k+P^{lj}_kP^{hi}_l-
P^{li}_kP^{hj}_l)\delta_h,
$$
$$
K(\delta_i,\delta_j)\partial^k=QQP^k_{ijh}\partial^h = (-R^k_{hij}
- R^0_{lij}Q^{lk}_h+S_{hil}P^{lk}_j-S_{hjl}P^{lk}_i)\partial^h,
$$
$$
K(\delta_i,\delta_j)\delta_k=QQQ_{ijk}^h\delta_h =(R^h_{kij}+
S_{ljk}P^{hl}_i-S_{lik}P^{hl}_j-R^0_{lij}P^{hl}_k)\delta_h,
$$
$$
K(\partial^i,\delta_j)\delta_k=PQQ^i_{jkh}\partial^h =(\partial^i
S_{hjk}+S_{ljk}Q^{il}_h-S_{hjl}P^{li}_k)\partial^h,
$$
$$
K(\partial^i, \delta_j)\partial^k=PQP_j^{ikh}\delta_h =(\partial^i
P^{hk}_j+P^{hi}_lP^{lk}_j-Q^{ik}_lP^{hl}_j)\delta_h.
$$

The explicit expressions of these components are obtained after
some quite long and hard computations, made by using the package
RICCI.

Next, the components of the Ricci tensor field are obtained as
traces of $K$
$$
RicPP^{jk}=Ric(\partial^j,\partial^k)=PPP^{hjk}_h-PQP^{jkh}_h,
$$
$$
RicQQ_{jk}=Ric(\delta_j,\delta_k)=QQQ_{hjk}^h+PQQ^h_{jkh},
$$
$$
Ric(\partial^j,\delta_k)=Ric(\delta_k,\partial^j)=0.
$$

\section{The cotangent bundle $T^*M$ as a K\"ahler Einstein
manifold}

From the explicit expressions of the components of the Ricci
tensor field on $T^*M$ one obtains the common Einstein factor
$Ef$, appearing in the condition for the K\"ahlerian manifold
$(T^*M,G,J)$ to be an Einstein space

$$
Ef= -n \frac{a_1^2a_1^\prime\lambda - 2a_1c\lambda +
a_1^3\lambda^\prime + 2a_1^\prime c\lambda t - 2a_1c\lambda^\prime
t}{2a_1\lambda^2(a_1 - 2a_1^\prime t)}\ -
$$

$$
-(a_1^2 - 2ct)(a_1a_1^\prime\lambda^2 + a_1^2\lambda\lambda^\prime
- a_1^{\prime 2}\lambda^2t + a_1a_1^{\prime\prime}\lambda^2t -
a_1^2\lambda^{\prime 2}t +
       a_1^2\lambda\lambda^{\prime\prime}t - 2a_1^{\prime 2}
       \lambda\lambda^\prime t^2 +
$$
$$
       2a_1a_1^{\prime\prime}\lambda\lambda^\prime t^2 +
       2a_1a_1^\prime\lambda^{\prime 2}t^2 -
       2a_1a_1^\prime \lambda \lambda^{\prime \prime} t^2)/
       (a_1\lambda^2(a_1 - 2a_1^\prime t)^2(\lambda + 2\lambda^\prime
   t))
$$

Next we consider the differences
$$
DiffQQ_{jk}=RicQQ_{jk}-Ef\ G^{(1)}_{jk}=
$$
$$
=\frac{a_1^2-2ct}{2a_1^2\lambda^2(a_1-2a_1^\prime
t)^4(\lambda+2\lambda^\prime t)^2}\gamma p_jp_k,
$$

$$
DiffPP^{jk}=RicPP^{jk}-Ef\ G_{(2)}^{jk}=
$$
$$
=\frac{1}{2a_1^2\lambda^2(a_1-2a_1^\prime
t)^2(a_1^2-2ct)(\lambda+2\lambda^\prime t)^2}\gamma g^{0j}g^{0k},
$$ \vskip2mm
\noindent where  $\gamma = nC_n+\beta$ and $C_n, \beta$ are
expressions involving $a_1,\lambda$ and their derivatives up to
third order. The condition for $(T^*M,G,J)$ to be K\"ahler
Einstein is given by $DiffQQ_{jk}=0, DiffPP^{jk}=0$ or,
equivalently, $\gamma=0$. If we ask for K\"ahler Einstein
structures on $T^*M$ to be independent of the dimension $n$ of
$M$, then we must have $C_n=0, \beta=0$. The coefficient of $n$ in
the expression of $\gamma$ is
$$
C_n = -(a_1 - 2a_1^\prime t)(a_1^2 - 2ct)(\lambda +
2\lambda^\prime t)^2 (2a_1a_1^{\prime 2}\lambda^2 +
a_1^2a_1^{\prime\prime}\lambda^2 +
2a_1^2a_1^\prime\lambda\lambda^\prime -
$$
$$
- 2a_1^3\lambda^{\prime 2}+ a_1^3\lambda\lambda^{\prime\prime} -
2a_1^{\prime 3}\lambda^2 t -
    2a_1a_1^{\prime 2}\lambda\lambda^\prime t +
    2a_1^2a_1^{\prime\prime}\lambda\lambda^\prime t +
    4a_1^2a_1^\prime\lambda^{\prime 2} t -
    2a_1^2a_1^\prime\lambda\lambda^{\prime\prime}t).
$$

Excluding the cases for which we have singularities, we can obtain
from the condition $C_n = 0$ the expression
\begin{equation}
\left\{
\begin{array}{c}
 a_1^{\prime\prime}= -(2a_1a_1^{\prime 2}\lambda^2 +
2a_1^2a_1^\prime\lambda\lambda^\prime - 2a_1^3\lambda^{\prime 2} +
a_1^3\lambda\lambda^{\prime\prime} - 2a_1^{\prime 3}\lambda^2 t- \\
-2a_1 a_1^{\prime 2}\lambda\lambda^\prime t +
4a_1^2a_1^\prime\lambda^{\prime 2}t -
2a_1^2a_1^\prime\lambda\lambda^{\prime\prime}t)/(a_1^2\lambda^2 +
2a_1^2\lambda\lambda^\prime t).
\end{array}
\right.
\end{equation}

Differentiating the expression of $a_1^{\prime\prime}$ with
respect to $t$, then replacing $a_1^{\prime\prime}$ from (15), we
get a quite complicate expression for the derivative of third
order $a_1^{(3)}$
$$
a_1^{(3)}= (12a_1^2a_1^{\prime 3}\lambda^4 - 24ta_1a_1^{\prime 4}
\lambda^4 + 12t^2a_1^{\prime 5}\lambda^4 + 18a_1^3a_1^{\prime
2}\lambda^3\lambda^\prime - 30ta_1^2a_1^{\prime
3}\lambda^3\lambda^\prime +
$$
$$
    +12t^2a_1a_1^{\prime 4}\lambda^3\lambda^\prime +
    18ta_1^3a_1^{\prime 2}\lambda^2\lambda^{\prime 2} -
    24t^2a_1^2a_1^{\prime 3}\lambda^2\lambda^{\prime 2} -
    12a_1^5\lambda\lambda^{\prime 3} +
    36ta_1^4a_1^\prime\lambda\lambda^{\prime 3} -
    $$
    $$
    -24t^2a_1^3a_1^{\prime 2}\lambda\lambda^{\prime 3} -
    12ta_1^5\lambda^{\prime 4} + 24t^2a_1^4a_1^\prime\lambda^{\prime 4} +
    3a_1^4a_1^\prime\lambda^3\lambda^{\prime\prime} -
    12ta_1^3a_1^{\prime 2}\lambda^3\lambda^{\prime\prime} +
    $$
    $$
    +12t^2a_1^2a_1^{\prime 3}\lambda^3\lambda^{\prime\prime}
    +9a_1^5\lambda^2\lambda^\prime\lambda^{\prime\prime} -
    24ta_1^4a_1^\prime\lambda^2\lambda^\prime\lambda^{\prime\prime} +
    12t^2a_1^3a_1^{\prime 2}\lambda^2\lambda^\prime\lambda^{\prime\prime} +
    $$
    $$
    +12ta_1^5\lambda\lambda^{\prime 2}\lambda^{\prime\prime} -
    24t^2a_1^4a_1^\prime\lambda\lambda^{\prime 2}\lambda^{\prime\prime}
    -a_1^5\lambda^3\lambda^{(3)} +
    2ta_1^4a_1^\prime\lambda^3\lambda^{(3)} - 2ta_1^5\lambda^2\lambda^\prime\lambda^{(3)} +
    $$
    $$
    +4t^2a_1^4a_1^\prime\lambda^2\lambda^\prime\lambda^{(3)})/
    (a_1^4(\lambda^2 + 2t\lambda\lambda^\prime)^2).
$$

Next, replacing these expressions of $a_1^{\prime\prime}, \
a_1^{(3)}$ in the condition $\beta=0$, we get the following
interesting relation

\begin{equation}
\left\{
\begin{array}{c}
\lambda (a_1-2ta_1^\prime)^3(a_1^2-2ct) (a_1\lambda-ta_1^\prime
\lambda+t
a_1\lambda^\prime)(a_1^\prime\lambda+\\
+ a_1\lambda^\prime)(a_1^2a_1^\prime\lambda+2c
a_1\lambda+a_1^3\lambda^\prime -2 c t a_1^\prime\lambda+2 c t a_1
\lambda^\prime)=0.
\end{array}
\right.
\end{equation}

The vanishing of the factors $\lambda , a_1-2ta_1^\prime$ and
$a_1^2-2ct$ will be not considered since the corresponding
situations lead to singularities. Thus we have the following three
essential cases \vskip2mm

1) The first and most interesting situation which will be studied
is that when the last factor in (16) vanishes. From the
corresponding relation one gets
\begin{equation}
a_1^\prime=-\frac{2ca_1\lambda+a_1^3\lambda^\prime+2ct
a_1\lambda^\prime}{\lambda(a_1^2-2ct)}.
\end{equation}

Differentiating $a_1^\prime$ with respect to $t$ and replacing
$a_1^\prime$ from (17) in the obtained result, one gets the same
expression for $a_1^{\prime\prime}$ as that obtained from (15),
after the replacing of $a_1^\prime$ from (17). Next, computing
$a_1^{(3)}$ and replacing again $a_1^\prime$ from (17) one gets
$DiffQQ_{jk}=0$ and $DiffPP^{jk}=0$. Thus if the relation (17) is
fulfilled, one obtains that $(T^*M,G,J)$ is K\"ahler Einstein.
Next one obtains the expresion
$$
Ef = 2c(n+1)\frac{a_1}{\lambda(a^2+2ct)},
$$
of the Einstein factor which must be a constant. We shall take
$Ef= \frac{k(n+1)}{2}$, where $k$ is a constant. It follows that
we can express $\lambda$ as a function of $a_1$ (although the
above computations could suggest to express $a_1$ as a function of
$\lambda$)

\begin{equation}
\lambda =\frac{4c}{k}\frac{a_1}{a_1^2+2ct}.
\end{equation}

Differentiating (18) with respect to $t$ it follows that (17) is
identically fulfilled. Hence the expression (18) of $\lambda$ is
obtained from a prime integral of (17).

{\bf Remark}. The same result is obtained if we express from the
equation $C_n=0$ the derivative $\lambda^{\prime\prime}$ as a
function of $\lambda, \lambda^\prime, a_1, a_1^\prime ,
a_1^{\prime\prime}$. \vskip2mm

Recall that the K\"ahler manifold $(T^*M,G,J)$ has constant
holomorphic sectional curvature $k$ if its curvature tensor field
$K$ can be expressed  by the relation
$$
K(X,Y)Z=\frac{k}{4}(G(Z,Y)X-G(Z,X)Y+
$$
$$
+G(Z,JY)JX-G(Z,JX)JY+ 2G(X,JY)JZ),
$$
where $X,Y,Z$ are vector fields on $T^*M$.

We shall use an adapted local frame $(\partial^1,...,\partial^n,~
\delta_1,..., \delta_n)$ in order to obtain the expressions for
the components of $K$ in the case where $(T^*M,G,J)$ has constant
holomorphic curvature. Introduce the following $M$-tensor fields
$$
J^{(1)}_{ij}=a_1g_{ij}+b_1p_ip_j,\ J_{(2)}^{kl}=a_2g^{kl}+
b_2g^{0k}g^{0l}.
$$
Remark that, up to a sign, the $M$-tensor fields $J^{(1)}_{ij},
J_{(2)}^{kl}$ are the components of the tensor field $J$, defining
the integrable almost complex structure on $T^*M$. Next we have
$$
K(\delta_i,\delta_j)\delta_k=\frac{k}{4}(G^{(1)}_{jk}\delta^h_i-
G^{(1)}_{ik}\delta^h_j)\delta_h, \ \
K(\partial^i,\partial^j)\partial^k=\frac{k}{4}(G_{(2)}^{kj}\delta_h^i-
G_{(2)}^{ki}\delta_h^j)\partial^h,
$$
$$
K(\delta_i,\delta_j)\partial^k=\frac{k}{4}(J^{(1)}_{ih}J^{(1)}_{jl}-
J^{(1)}_{il}J^{(1)}_{jh})G_{(2)}^{kl}\partial^h,
$$
$$
K(\partial^i,\partial^j)\delta_k=\frac{k}{4}(J_{(2)}^{ih}J_{(2)}^{jl}-
J_{(2)}^{il}J_{(2)}^{jh})G^{(1)}_{kl}\delta_h,
$$
$$
K(\partial^i,\delta_j)\partial^k=\frac{k}{4}(-J^{(1)}_{jl}J_{(2)}^{ih}
G_{(2)}^{kl}-
G_{(2)}^{ki}\delta^h_j-2J^{(1)}_{jl}J_{(2)}^{kh}G_{(2)}^{il})\delta_h,
$$
$$
K(\partial^i,\delta_j)\delta_k=\frac{k}{4}(G^{(1)}_{kj}\delta_h^i+
G^{(1)}_{kl}J_{(2)}^{il}J^{(1)}_{jh}+2
G_{(2)}^{il}J^{(1)}_{jl}J^{(1)}_{kh})\partial^h
$$

In our case, i.e. when $a_1,\lambda$ are related by (18), one
obtains that the components of $K$ are given just by the above
relations, hence the K\"ahler Einstein manifold $(T^*M,G,J)$ has
constant holomorphic curvature $k$. Hence we may state the
following result.

{\bf Theorem 9.} \it Assume that the Riemannian manifold $(M,g)$
has constant sectional curvature $c$ and consider the natural
integrable almost complex structure $J$ defined on its cotangent
bundle $T^*M$ by $(3)$, where the coefficients $a_1,a_2, b_1,b_2$
are related by $(4)$ and $(7)$. There exists a family of K\"ahler
Einstein structures $(G,J)$  defined by (9), on $T^*M$, where the
coefficients $c_1,c_2,d_1,d_2$ are expressed by $(11),\ (13)$ and
the factors $\lambda , \ \mu$ are given by $\mu = \lambda^\prime$
and by $(18)$. Moreover, the obtained K\"ahler Einstein structure
has constant holomorphic sectional curvature $k$. \rm

\bf Remark. \rm The parameter $a_1$ is not quite arbitrary. In
fact, the following conditions must be fulfilled

$$
a_1>0,\ a_1+2tb_1 = \frac{a_1a_1^\prime-c}{a_1-2ta_1^\prime}>0,
$$
$$
\lambda =\frac{4c}{k}\frac{a_1}{a_1^2+2ct}>0,\ \lambda
+2t\lambda^\prime =
\frac{4c}{k}\frac{(a_1-2ta_1^\prime)(a_1^2-2ct)}{(a_1^2+2ct)^2}>0.
$$

{\bf Example.} Assume $c>0$ and consider the function
$a_1=B+\sqrt{B^2+2ct}$, where $B$ is a positive constant. We have
$a_1^\prime =\frac{c}{\sqrt{B^2+2ct}}$ and one checks easily that
all the conditions from the above remark are fulfilled. In the
case where $c<0$ one can consider the function
$a_1=B+\sqrt{B^2-2ct}$, where $B$ is a positive constant, and a
simple algebraic computation shows that $\lambda$ is a constant
and all the conditions from the above remark are fulfilled. In
fact, the case $\lambda =1$ has been considered, in the case of
the tangent bundle, in \cite{Oproiu4}, \cite{Oproiu5}.\vskip3mm

2) The next situation is obtained when
$$
a_1^\prime\lambda+a_1\lambda^\prime=0
$$

It follows that $a_1\lambda =k$, a constant (this constant is not
related to the constant used in the study of the first case). Then
we have
$$
\lambda=\frac{k}{a_1},\ \lambda^\prime=- k
\frac{a_1^\prime}{a_1^2},\ \lambda^{\prime\prime} =k \frac{2
a_1^{\prime 2}-a_1 a_1^{\prime\prime}}{a_1^3},
$$
$$
\lambda^{(3)}=k\frac{6 a_1a_1^\prime
a_1^{\prime\prime}-6a_1^{\prime 3}- a_1^2a_1^{(3)}}{a_1^4}
$$

With these values of $\lambda,\lambda^\prime,
\lambda^{\prime\prime}, \lambda^{(3)}$ one gets that the
conditions \break $DiffQQ_{jk}=0$ and $DiffPP^{jk}=0$ are
fulfilled identically.

If we study the property of $(T^*M, G,J)$ to have constant
holomorphic sectional curvature, we get that the components
$K(\delta_i,\delta_j)\delta_k, \ K(\delta_i,\delta_j)\partial^k$,
$\ K(\partial^i,\partial^j)\delta_k, \
K(\partial^i,\partial^j)\partial^k$ can be expressed just like in
the case 1. However, the last two components
$K(\partial^i,\delta_j)\partial^k, \
K(\partial^i,\delta_j)\delta_k$ are quite different from the
expression obtained in the case 1. Hence $(T^*M, G,J)$ cannot have
constant holomorphic sectional curvature. Then, we may state the
following result

{\bf Theorem 10.}  {\it Consider the K\"ahlerian structure $(G,J)$
on $T^*M$ obtained in theorem 7, depending on the essential
parameters $a_1,\lambda$. If $\lambda =\frac{k}{a_1}$, then the
manifold $(T^*M,G,J)$ is K\"ahler Einstein. It cannot have
constant holomorphic sectional curvature.}

{\bf Remark.} In this case, the following conditions must be
fulfilled
$$
a_1>0,\ a_1+2tb_1= \frac{a_1a_1^\prime-c}{a_1-2ta_1^\prime}>0,
$$
$$
k>0,\  \frac{a_1^2}{k}(\lambda+2t\lambda^\prime)=
a_1-2ta_1^\prime>0,
$$
hence the functions $a_1^2-2ct$ and $\frac{t}{a_1^2}$ must be
increasing. \vskip3mm

 3. The last case is obtained when
 $$
 a_1\lambda-ta_1^\prime \lambda+t a_1\lambda^\prime=0
 $$

 One sees easily that, in this case, one has
 $$
a_1=kt\lambda,
 $$
where $k$ is a constant (this constant is not related to the
constant used in the previous cases). One sees easily that
$a_1(0)=0$, thus this situation should be excluded. However, we
can study the properties of the K\"ahlerian structure $(G,J)$ on
the manifold $T_0^*M$ obtained from $T^*M$ by excluding the zero
section. Next one gets that the conditions $DiffQQ_{jk}=0$ and
$DiffPP^{jk}=0$ are fulfilled identically, so the K\"ahler
manifold $(T_0^*M,G,J)$ is Einstein.

{\bf Theorem 11.}  {\it Consider the K\"ahlerian structure $(G,J)$
on $T^*M$, obtained in theorem 7, depending on the essential
parameters $a_1,\lambda$. If  $a_1=kt\lambda$ then the manifold
$(T^*_0M,G,J)$ is K\"ahler Einstein.}

{\bf Remark.} The function $\lambda$ must fulfill the conditions
obtained from $a_1>0, a_1+2tb_1>0, \lambda>0,\
\lambda+2t\lambda^\prime>0$.

\vskip 1.5cm
\begin{minipage}{2.5in}
\begin{flushleft}
V.Oproiu\\
Faculty of Mathematics\\
University "Al.I.Cuza", Ia\c si \\
Rom\^ania.\\
e-mail: voproiu@uaic.ro
\end{flushleft}
\end{minipage}
\hfill
\begin{minipage}{2.5in}
\begin{flushleft}
D.D.Poro\c sniuc\\
Department of Mathematics\\
National College "M. Eminescu" \\
Boto\c sani, Rom\^ania.\\
e-mail: danielporosniuc@lme.ro \\
~~~~~~~~~~dporosniuc@yahoo.com
\end{flushleft}
\end{minipage}

\end{document}